\documentclass[ a4paper,11pt]{amsart}
\usepackage{amssymb,amsfonts,epsf,epsfig,here,
rotate,latexsym}

\def\cl{\centerline}
\def\bs{\vskip24pt}
\def\ms{\vskip12pt}
\def\ss{\vskip6pt}
\def\nt{\noindent}
\begin{document}
\title{Discussion on Benford's Law and Its Application}
\author{Li Zhipeng, Cong Lin and Wang Huajia}
\address{\noindent{\small
Hwa Chong Junior College\\
661 Bukit Timah Road\\
Singapore 269734}} \email{phy@scientist.com (Li Zhipeng)\\
ouranos-veritas@scientist.com (Cong Lin)}
\date{}
\maketitle

\begin{abstract}
The probability that a number in many naturally occurring tables
of numerical data has first significant digit (i.e., first
non-zero digit) $d$ is predicted by Benford's Law\ ${\rm
Prob}\,(d) = \log_{10} \left(1 + {\displaystyle{1\over
d}}\right),  d = 1, 2 \ldots, 9$. Illustrations of Benford's Law
from both theoretical and real-life sources on both science and
social science areas are shown in detail with some novel ideas and
generalizations developed solely by the authors of this paper.
Three tests, Chi-Square test, total variation distance, and
maximum deviations are adopted to examine the fitness of the
datasets to Benford's distribution. Finally, applications of
Benford's Law are summarized and explored to reveal the power of
this mathematical principle.
\end{abstract}
\bs \nt {\bf 1. \ \ Introduction} \ss The significant-digit law of
statistical folklore is the empirical observation that in many
naturally occurring tables of numerical data, the leading
significant digits are not uniformly distributed as might be
expected, but instead follow a particular logarithmic
distribution. Back in 1881, the astronomer and mathematician Simon
Newcomb published a 2-page article in the American Journal of
Mathematics describing his observation that books of logarithms in
the library were dirtier in the beginning and progressively
cleaner throughout [1]. He inferred that researchers, be them
mathematicians, biologists, sociologists as well as physicists,
were looking up numbers starting with 1 much more often than
numbers beginning with 2, and numbers with first digit 2 more
often than 3, and so on. This ingenious discovery led him to
conclude that the probability that a number has first significant
digit (i.e., first non-zero digit) $d$ is ${\rm Prob}\, (d) =
\log_{10} \left(1 + {\displaystyle{1\over d}}\right), \ \ d = 1, 2
\ldots 9$. In particular, his conjecture stated that the first
digit is 1 about 30.1\% of the time, and is 9 only about 4.6\% of
the time. That digits are not equally likely comes as somewhat of
a surprise, but to claim an exact law describing their
distribution is indeed striking. Passed unnoticed, the proposed
law was discovered again and supported by empirical evidence by
the physicist Benford who analyzed the frequencies of significant
digits from twenty different tabled including such diverse data as
surface areas of 335 rivers, specific heat of thousands of
chemical compounds, and square-root tables [2]. And this First
Digit Law is known as Benford's Law today. But in recognition of
Newcomb's discovery, we can call it Newcomb-Benford's Law. This
law applies to stock prices [3], number of hours billed to clients
[4] or income tax [5] as well as mathematical series [6]. And the
tremendous practical values of Benford's Law were neglected until
recently many mathematicians began to focus on the applications of
this amazing phenomenon such as the design of computers and
analysis of roundoff errors [7, 14, 15], as well as a
goodness-of-fit against Benford to detect fraud [8]. \ss In this
paper, we will be dealing with heuristic argument and
distributional property of the Significant-digit law in section 2;
checking how data from various sources fit Benford's Law in
section 3; discussing the application of Benford's law in section
4. \bs \nt{\bf  2. \ \ Mathematical Formulation} \ss \nt{\bf 2.1 \
\ Heuristic Argument} \ss Pietronero and his colleagues gave a
general explanation for the origin of the Benford's law in terms
of multiplicative processes in 2001 [3]. Here, the explanation was
amended slightly in such a way that it can be used to explain the
Benford's distribution not specifically to base 10. It stated that
many systems such as the stock market prices which is discussed
later do not follow the dynamical description by a Brownian
process:
$$N(t+1) = \xi + N(t),$$
but rather a multiplicative process:
$$N(t+1)=\xi N(t)$$
where $\xi$ is a stochastic variable. By a simple transformation, we get
$$\ln{N(t+1)} = \ln{\xi} + \ln{N(t)}.$$
If we consider $\ln \xi$ as the new stochastic variable, we
recover a Brownian dynamics in a logarithmic space; here we mean
that a random multiplicative process corresponds to a random
additive process in logarithmic space. This implies that as $t\to
\infty$, the distribution ${\rm Prob}\,(\ln N)$ approaches a
uniform distribution. By transforming back to the linear space we
have
$$\int \,{\rm Prob}\,(\ln N) d (\ln N) = \int C d(\ln N) = C \int {1\over N} dN,$$
where $C$ is a constant. \ss It should be noted that ${\rm Prob}\,
(N) = {\displaystyle{1\over N}}$ is not a proper probability
distribution, as it diverges or put in another way $\int_0^\infty
{\displaystyle{1\over N}}\,dN$ is undefined. However, the physical
laws and human conventions usually impose maximums and minimums.
The probability that the first significant digit of $N$ is $n$ in
base $b$ is given by the following expression:
$${\rm Prob} \, (n) = {\int_n^{n+1} {1\over N} dN\over \int_1^b {1\over N} dN} = {\ln {n+1\over n}\over \ln b} = \log_b \left(1 + {1\over n}\right),$$
for any integer $n$ which is less than $b$. We can review ${\rm
Prob}\,(n) = \log_b \left(1 + {\displaystyle{1\over n}}\right)$ as
a generalized expression of Benford's law to arbitrary base $b$.
\ms \nt{\bf 2.2 \ \ The Significant-Digit Law and Some
Consequences} \ss \nt {\bf 2.2.1 \ \ Significant-Digit Law} \ss
The Significant-Digit Law is
$${\rm Prob}\, (D_1 = d_1, \ldots, D_k = d_k) = \log_{10} \left[ 1 + \left( \sum_{i=1}^k \; d_k \times 10^{k-1}\right)^{-1}\right]$$
where $D_1, D_2,\ldots D_k$ are the first, second $\ldots$ $k$'th
digits respectively. \ss For example, ${\rm Prob} \, (D_1 = 1, D_2
= 2, D_3 = 9) = \log_{10} (1 + (129)^{-1}) \cong 0.00335$, which
means there is a probability of 0.00335 that the first three
significant digits are 129 in a sample of Benford's distribution
[11]. \ss Hill has proved in his papers: ``Scale-Invariance
implies Base-Invariance'', not vice versa [9]; ``Base-Invariance
implies Benford's Law [10]''; ``The logarithmic distribution is
the unique continuous base-invariant distribution [9]''. He has
explained the Central-limit-like Theorem for Significant Digit by
saying: ``Roughly speaking, this law says that if probability
distributions are selected at random and random samples are then
taken from each of these distributions in any way so that the
overall process is scale (or base) neutral, then the
significant-digit frequencies of the combined sample will converge
to the logarithmic distribution [11]''. \ms \nt{\bf 2.2.2 \ \
Distributional Properties of $D_k$'s} \ss Based on the
Significant-Digit Law, we analyze the statistics and
distributional properties of a Benford's Distribution. \ss \nt
{\bf 2.2.2.1 \ \ Mean and Variance for $D_k$} \ss In this
subsection, we compute the numerical values of the means and
variances of $D_k$'s using these expressions: \begin{equation}
E(D_k) = \sum_1^9 \; n \,{\rm Prob}\, (D_k = n) \nonumber
\end{equation} \begin{equation}{\rm Var}\, (D_k) = \sum_{n=1}^9 \; n^2 \,{\rm
Prob}\, (D_k = n) - E (D_k)^2,\nonumber \end{equation} and
tabulate them below.

\ms
{\small \cl{Table 1. Mean and variance of $D_k$ for $k=1$ to 7}
\hspace{1.25in}\begin{tabular}{|c|ll|ll|}\hline
$k$ & \multicolumn{2}{c}{$E(D_k)$} & \multicolumn{2}{c|}{\mbox{Var $(D_k)$}} \\ \hline
1 &3.44023696712 &(5)   &6.0565126313757 &(6.67)\\ \hline
2 &4.18738970693 &(4.5) &8.2537786232732 &(8.25)\\ \hline
3 &4.46776565097 &(4.5) &8.2500943647286 &(8.25)\\ \hline
4 &4.49677537552 &(4.5) &8.2500009523513 &(8.25)\\ \hline
5 &4.49967753636 &(4.5) &8.2500000095245 &(8.25)\\ \hline
6 &4.49996775363 &(4.5) &8.2500000000953 &(8.25)\\ \hline
7 &4.49999677536 &(4.5) &8.2500000000016 &(8.25)\\ \hline
\end{tabular}
}

\ms This shows that by the Significant-digit law the mean of $D_k$
is approaching 4.5 which is the mean if the distribution were
uniform and the variance of $D_k$ is approaching 8.25 which is the
variance if the distribution were uniform. \ms \nt{\bf 2.2.2.2 \ \
Histogram of $D_k$ for $k$}

\raisebox{5ex}{\includegraphics[width=1.808in,height=1.848in]{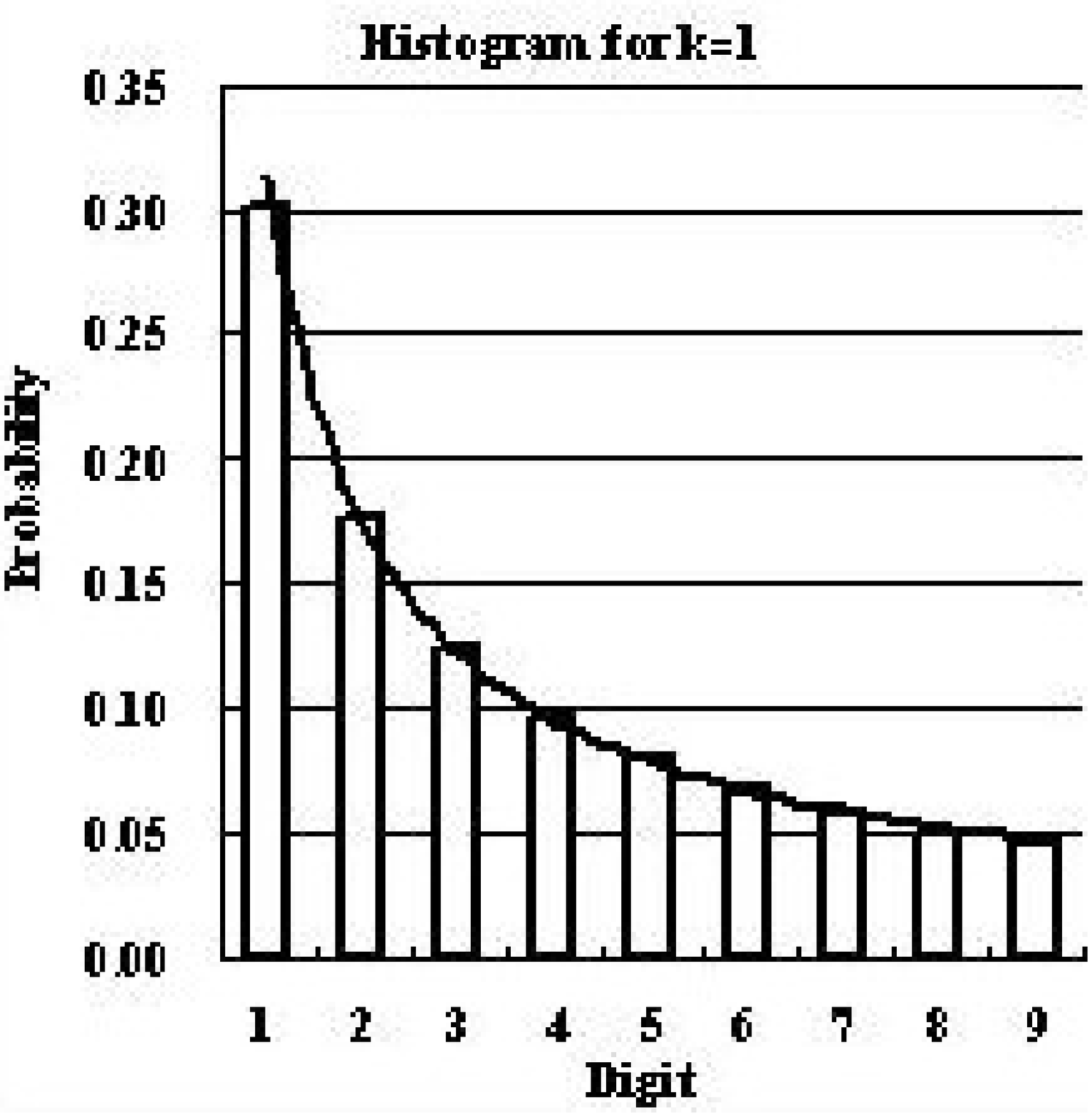}}
\hspace{0in}\includegraphics[width=3.752in,height=2.864in]{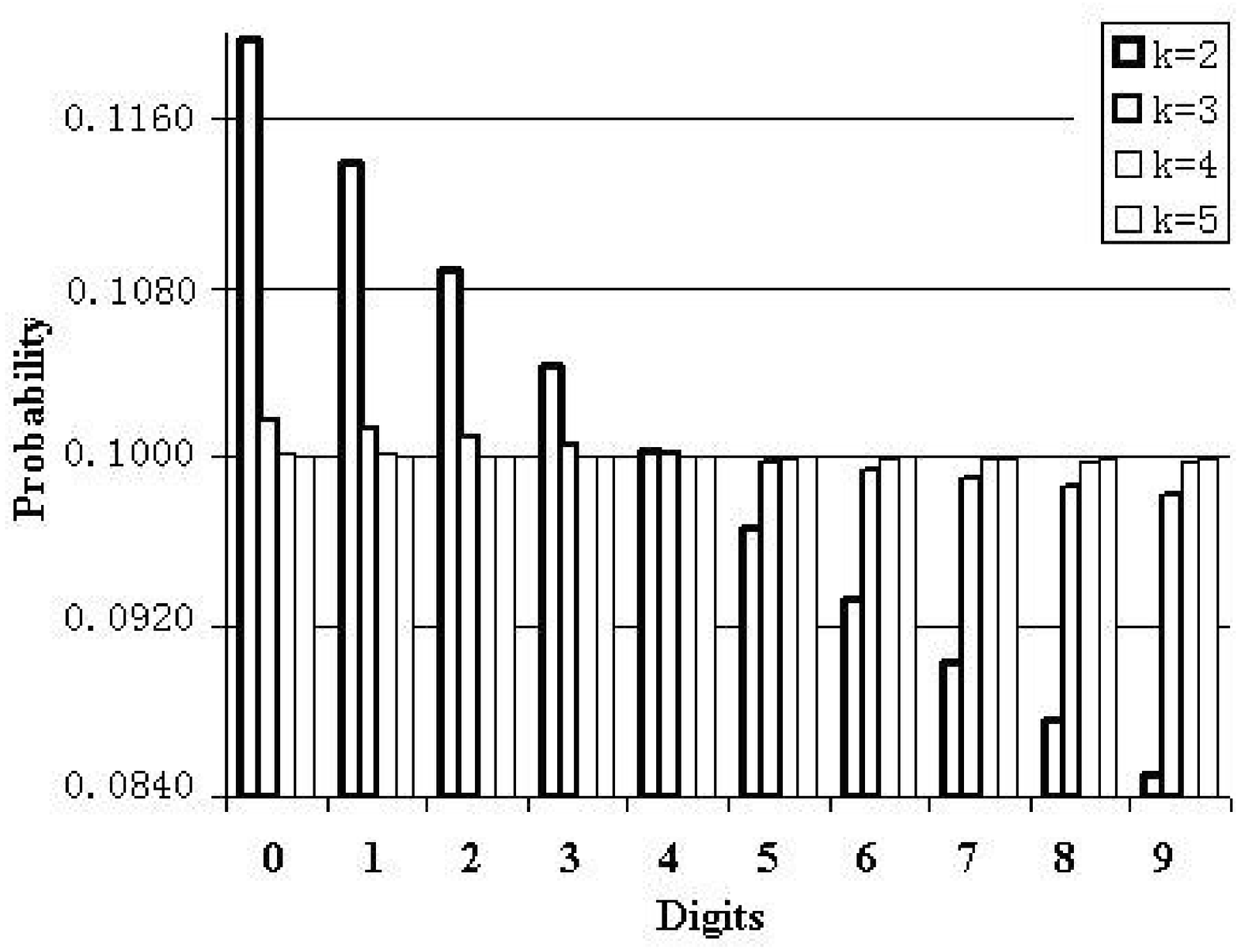}

\cl{\small Fig 1. Histograms of $D_k$ for $k=1$ \hspace{0.5in}Fig 2. Histograms of $D_k$ for $k=2,3,4,5$}
\ss
With the histogram, we illustrate that the distribution of the nth significant digit approaches the uniform distribution rapidly as $n\to\infty$.

\ms
\nt{\bf 2.2.2.3 \ \ Total variation distance of $D_k$ with respect to a uniform distribution}

\begin{minipage}[b]{3.5in} By calculating the total variation distance $d(D_1, U)$ of $D_k$ from a uniform distribution on $\{0, 1, 2 \ldots 9\}$ using these expressions:
$$d(D_1, U) = {1\over 2} \sum_{n=1}^9 \left| {\rm Prob} \, (D_1 = n) - {1\over 9} \right| \, {\rm for}\, k = 1$$
$$d(D_k, U) = {1\over 2} \sum_{n=0}^9 \left|{\rm Prob} \, (D_k = n) - {1\over 10} \right| \, {\rm for}\, k \not= 1.$$
As shown in the table, we observe that the $D_k$ is converging to
the uniform distribution geometrically. \vspace{-0.6in}
\end{minipage}\hspace{0.4in}
\begin{minipage}[b]{1.5in}{\small Table 2. Total Variation
Distance of $D_k$ from a
uniform distribution on
$\{0, 1, 2, \ldots, 9\}$}.
\ss
{\small \begin{tabular}{|l|l|} \hline
$k$ &$d(D_k,U)$\\ \hline
1 &0.26872666 \\ \hline
2 &0.04702863\\ \hline
3 &0.00488356\\ \hline
4 &0.00048858\\ \hline
5 &0.00004886\\ \hline
6 &0.00000489\\ \hline
7 &0.00000049\\ \hline
$k \to \infty$ &$d(D_k, U)\to 0$\\ \hline
\end{tabular}
}\end{minipage}

\ms
\nt{\bf 2.2.2.4 \ \ Correlation coefficient of $D_i$ and $D_j$}
\ss
Here, by calculating the correlation coefficient of $D_i$ and $D_j$,
$${\displaystyle \rho_{_{_{D_i D_j}}} = {{\rm Cov}\, (D_i, D_j) \over \sqrt{{\rm Var}\, (D_i) \, {\rm Var}\, (D_j)}}}$$

for $0 < i<j$. we intend to investigate the dependence of one digit on another digit. We conclude that the dependence among significant digits decreases as the distance $(j-i)$ increases.
\ss
{\small \cl{Table 3.Correlation coefficient of $D_i$ and $D_j$}
\hspace{1in}\begin{tabular}{|l|c|c|c|c|c|} \hline
\begin{minipage}[t]{0.3in}
\begin{picture}(10,10)(0,0)
\put (-5,10){\line(2,-1){30}}
\put (2.5,0.78){\makebox(0,0){\small $i$}}
\put (22,4.5){\makebox(0,0){\small $j$}}
\end{picture}
\end{minipage} &2 &3 &4 &5\\ \hline
1 &0.0560563 &0.0059126 &0.0005916 &0.0000591\\ \hline
2 & &0.0020566 &0.0002059 &0.0000205 \\ \hline
3 &&&0.0000228 &0.0000022\\ \hline
4 &&&&0.0000002\\ \hline
\end{tabular}
}

\ms
\nt{\bf 3. \ \ Illustrations}
\ss
The Chi-Square test formula to check the sample's goodness of fit to Benford's distribution is given as below,
$$\chi^2 (8^\circ ) = \sum_{n=1}^9 {(\log \left(1 + {1\over n}\right) - \,{\rm Prob}\, (D_1 = n))^2\over \log \left(1 + {1\over n}\right)} \times S,$$
where $S$ denotes the sample size. The critical values for 8
degree of freedom at 5\% and 1\% level of significance are
respectively 15.51 and 20.09. \ss Another test, which we have used
to analyze the sample, is the total variation distance (denoted as
$d_1$ here) of the sample from the Benford's distribution on $\{1,
2, \ldots, 9\}$:
$$d_1 = {1\over 2} \sum_{n=1}^0 \left| {\rm Prob}\, (D_1 = n) - \log (1 + {1\over n})\right|.$$
The maximum of deviations from a uniform distribution for each digit could also be considered as a test to
check the goodness of fit:
$$d_{\max} = {\mathop{\max}\limits_{1\le n\le 9}} \left\{\left| {\rm Prob}\, (D_1 = n) =
\log \left(1 + {1\over n}\right)\right|\right\}.$$
\ms
\nt{\bf 3.1 \ \ Physical Constants}
\ss
Many literatures on Benford's law cite the table of physical constants as an example to illustrate
Benford's law [9, 10, 11]. It was only Burke and his colleagues who actually attempted to check whether
the physical constants would match Benford's law [12].

\ms
\begin{minipage}{2.4in}
However, they only chose the constants from an introductory
physics text. The sample size is too small to be significant. Here
the 183 constants from the 1998 Committee On Data for Science and
Technology recommended complete listing of the fundamental
physical constants (http://physics.nist.gov/constants) were
analyzed for the first time.
\end{minipage}\hspace{0.2in}
\begin{minipage}{2.9in}
\includegraphics[width=2.853in,height=1.665in]{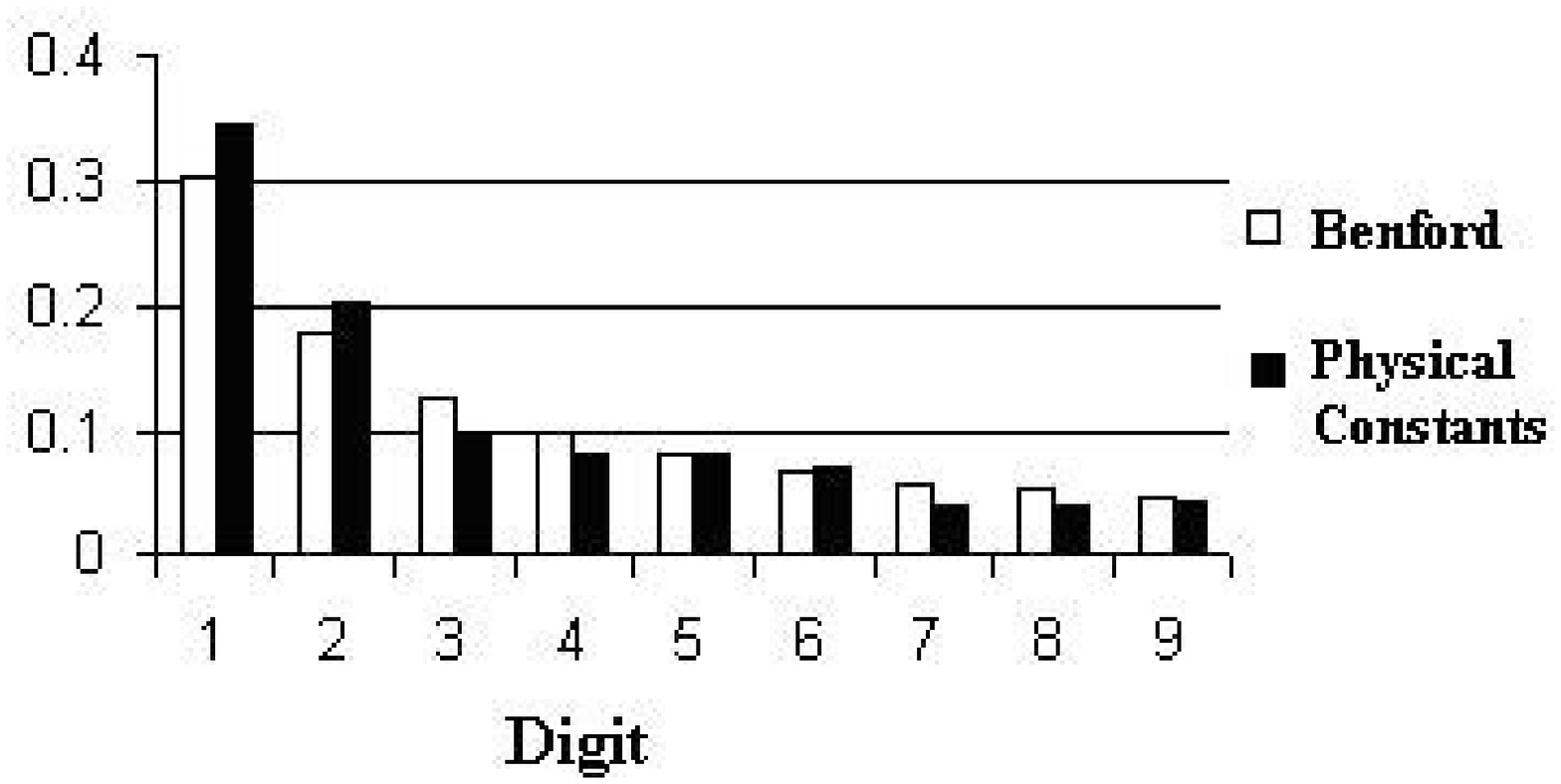}
\hspace*{-0.1in}{\small Fig 3. Histogram on fitness of physical constant}
\end{minipage}

\ms

{\small \cl{Table 4. First Significant Digit Test for the physical constants.}
\hspace{0.1in}\begin{tabular}{|l|l|l|l|l|} \hline
Digits &Counts in  &Predicted by  &Frequency
Observed &The Difference\\
&Sample Frequency &Benford's Law &in Sample &\\ \hline
1 &63 &0.3010 &0.3443 &0.0432 \\ \hline
2 &37 &0.1761 &0.2022 &0.0261\\ \hline
3 &18 &0.1249 &0.0984 &-0.0266\\ \hline
4 &15 &0.0969 &0.0820 &-0.0149\\ \hline
5 &15 &0.0792 &0.0820 &0.0028\\ \hline
6 &13 &0.0669 &0.0710 &0.0041\\ \hline
7 & 7 &0.0580 &0.0383 &-0.0197\\ \hline
8 &7 &0.0512 &0.0383 &-0.0129\\ \hline
9 &8 &0.0458 &0.0437 &-0.0020\\ \hline
\multicolumn{3}{|l}{Sample Size: 183} &\multicolumn{2}{l|}{Chi-Square Test: 5.206}\\ \hline
\multicolumn{3}{|l}{Total variation distance: 0.0762} &\multicolumn{2}{l|}{ Maximum of deviations: 0.0432}\\ \hline
\end{tabular}}

\bigskip
\nt{\bf 3.2 \ \ Stock Prices \& One-day Returns on Stock Index}
\ss Interesting results on the frequency of first significant
digit of one-day returns on stock index were found by Ley [13].
However, investigations on stock prices were not done. According
to Pietronero, due to the multiplicative process, a stock's prices
over a long period of time should conform to Benford's Law. Here,
we look at all stock's prices on certain days rather than a single
stock's prices on a series of days. We collected data on local
stock market SGX main board from http://www.sgx.com. \ms

{\small \cl{Table 5. Analysis of Stock Prices}
\hspace{0.5in}\begin{tabular}{|c|l|l|l|l|l|}\hline
    &Sample Source &Sample Size &$\chi^2$ &$d_1$ &$d_{\max}$ \\ \hline
$A$ &One-trading-day  &548 &2.23 &0.0277 &0.01734\\
&(13/10/2003)&&&&\\ \hline
$B$ &20-trading-Day  &11,015 &45.5 &0.0255 &0.01501 \\
&(10/10/2003- 7/11/2003)&&&&\\ \hline
$C$ &32-trading-Day &17,214 &79.3 &0.0284 &0.01834\\
&(10/10/2003- 26/11/2003) &&&& \\ \hline
\end{tabular}}

\bigskip
The observed frequencies roughly agree with the theoretical
frequency predicted by Benford's law, indicted by the relatively
small values of the total variation distance (0.0277, 0.0255 and
0.0284) of the sample and the Benford's distribution. However, if
we performed the usual chi-square test of goodness of fit, we
would reject the null hypothesis for datasets $B$ and $C$, as the
chi-square values are much greater than 15.51 for the critical
chi-square test value of 8 degrees of freedom at 5\% significant
level or 20.09 at 1\% significant level. It is because of the
large number of observations - that is, the classical acceptance
region shrinks with the sample size, given a significant level. On
the contrary, the Benford's hypothesis would not be rejected for
dataset $A$. \ms
\begin{minipage}{2.2in}
One interesting thing, which in three sets of data $d_{\max}$
always occurs at $n = 2$, attracts our attention. By plotting the
histograms of three datasets against a Benford's distribution, it
is found that the distributions of stock price data are in quite
agreement among themselves but a little different from Benford's
distribution.
\end{minipage}
\hspace{0.2in}
\begin{minipage}{3.24in}
\includegraphics[width=3.24in,height=2.28in]{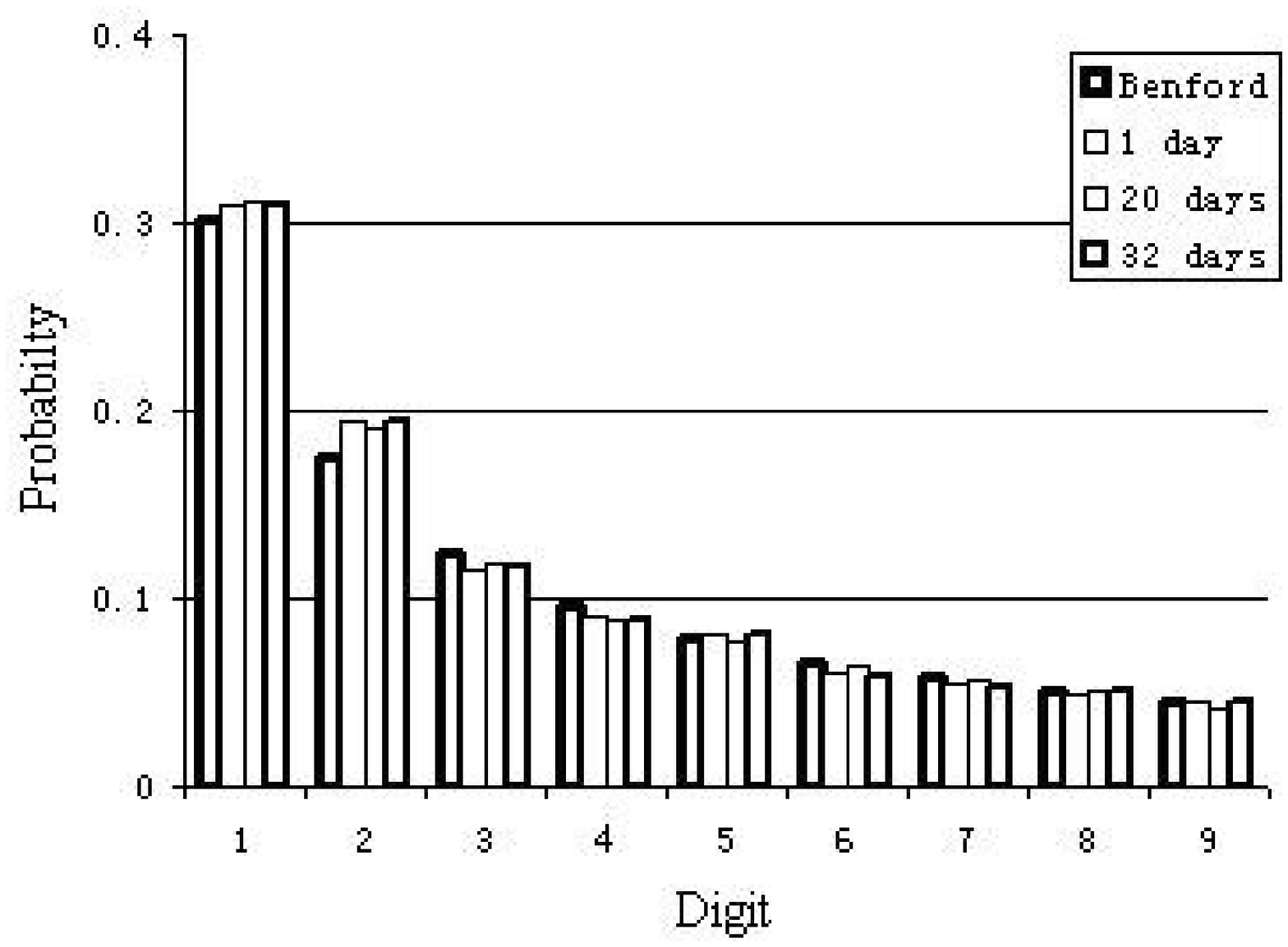}
\hspace{0.3in}{\small Fig 4. Histogram on fitness of stock prices}
\end{minipage}

We attribute this to the fact that a simple multiplicative process
cannot precisely reproduce the complicated stock process thus the
distribution of stock prices data may follow a curve that is
steeper than the log curve of Benford's Law.

\ms \nt{\bf 3.3 \ \ Some Mathematical Series} \ss \nt{\bf 3.3.1 \
\ Fibonacci Series} \ss Here we provide a simple explanation for
the conformance of Fibonacci series to Benford's law. Fibonacci
series is characterized by the recursion relation: $a_{n+2} =
a_{n+1} + a_n$.  This formula does not provide the explicit
information of a general term from the series that may suggest the
origin of Fibonacci series' conformance to Benford's law. However,
it could be easily derived that say for a series with $a_1 = 1$,
$a_2 = 2$:
$$a_n = \left({1 + \sqrt{5}\over 2}\right)^{n-1} \left({5+3\sqrt{5} \over 10}\right) + \left({1 - \sqrt{5}\over 2}\right)^{n-1} \left({5 - 3\sqrt{5}\over 10}\right).$$
Moreover, the magnitude of the second term is
${\left|\left({1-\sqrt{5}\over 2}\right)^{n-1}
\left({5-3\sqrt{5}\over 10}\right)\right|}$ which is always less
than 1 and approaches zero when $n$ gets larger and larger. Thus
to analyze the first significant digit of $a_n$, we only need to
analyze: $$a_n^\prime = \displaystyle{\left({1+\sqrt{5}\over
2}\right)^{n-1} \left({5 + 3\sqrt{5}\over 10}\right)}$$ instead.
Now we have ${\ln a_n^\prime = (n - 1) \ln \left({1 +
\sqrt{5}\over 2}\right) + \ln \left({5 + 3\sqrt{5}\over
10}\right).}$ When the sample size is large, distribution on $n$
could be considered as uniform, with $\ln \left({1 + \sqrt{5}\over
2}\right)$ and $\ln \left({5 + 3\sqrt{5}\over 10}\right)$ being
constants, $\{\ln  a_n^\prime\}$ will be uniform. It is then
followed that the $\{a_n^\prime\}$ and therefore $\{a_n\}$ will be
uniformly distributed in logarithmic space and will conform to
Benford's law. \ms \nt{\bf  3.3.2 \ \ The Prime-number Series} \ss
The prime number series is rather uniform below 100000, with the
probability of each possible first significant digit being between
12.5\% and 10.4\%. Moreover, using the upper and lower bounds of
function pi from the prime number theorem, it can be shown that
the prime number sequence approximates a uniform distribution. \ss
\nt{\bf 3.3.3 \ \ Sequence $\alpha^n$ in base $b$, where $\alpha
\in R$, $n \in N$} \ss The origin of Benford's distribution in
$\alpha^n$ is often attributed to its scale invariance which means
that any power law relation is scale invariant:
$${\rm Prob} \, (\lambda N) = f (\lambda) \,{\rm Prob}\, (N).$$
An even simpler argument can be done. For $\alpha^n$ in base $b$, where $\alpha$ is a constant, the sufficient condition for which the first significant digit of $\alpha^n$ is $d$ is $d \times b^k \le \alpha^n < (d+1)\times b^k$ (that is $\log_b d + k \le n \log_b \alpha < \log b (d+1 + k)$ for an integer $k$. Thus for each $k$, the probability for the first significant digit of $\alpha^n$ being $d$ is
$${(\log_b (d+1) + k) - (\log_b d + k)\over (k+1)-k} = \log_n \left({d+1\over d}\right).$$
\ms \nt{\bf 3.3.4 \ \ Factorial } \ss For the first 160 factorial
numerical values, the chi-square test does not reject the
Benford's distribution; however the total variation distance is
too big for it to be considered as fit. \ss \nt{\bf  3.3.5 \ \
Sequence of Power, i.e. $n^2, n^3, n^4, n^5, \ldots$} \ss For
$n^k$ series, as the constant $k$ increases, the distribution
becomes closer to the Benford's law. All values of the three tests
get smaller as $k$ increases, as shown in the analysis of $n^2,
n^5, n^{20}$ and $n^{50}$  for $n = 1$ to 30,000. \ss \nt {\bf
3.3.6 \ \ Numbers in Pascal Triangle} \ss We suspect that the
numbers in Pascal Triangle obeys Benford's law, because it relates
to Fibonacci Series closely while in another perspective it could
be considered as a mixture of sequences of power. However in this
paper, no data on Pascal triangle were analyzed. \ms

{\small \cl{Table 6. Analysis of some mathematical series}
\hspace{0in}\begin{tabular}{|l|l|l|l|l|} \hline
Sample & Sample Size &$\chi^2$ &$d_1$ &$d_{\max}$\\ \hline
Fibonacci Series (7 series) &10, 317 &0.0125 &$3.8 \times 10^{-4}$ &$1.7 \times 10^{-4}$\\ \hline
Prime-number below 1,000 &168 &45.0 &0.2271 & 0.1522 \\ \hline
Prime-number below 100,000 &9761 &3247 &0.4905 &0.1761 \\ \hline
$1.007^n$, Specifically for $n=1$ to 30,000  &30,000 &0.410 &$1.2 \times 10^{-3}$ &$5.9 \times 10^{-4}$ \\ \hline
$1.007^n$, Specifically for $n=1$ to 65,028 &65,028 &0.0329 &$2.5 \times 10^{-4}$ &$1,2 \times 10^{-4}$\\ \hline
Factorial of 1 to 100 & 100& 6.95& 0.0651&0.04885\\ \hline
Factorial of 1 to 130 & 130& 8.97& 0.0871& 0.03492\\ \hline
Factorial of 1 to 160 & 160& 10.10& 0.0834& 0.03635\\ \hline
$n^2$, Specifically for $n=1$ to 30,000 & 30,000&$3.16\times 10^3$ &0.1409 & 0.09900\\ \hline
$n^5$, Specifically for $n=1$ to 30,000  & 30,000& $2.76 \times 10^2$& 0.0394& 0.03640\\ \hline
$n^{20}$, Specifically for $n=1$ to 30,000 & 30,000& 20.8& 0.0112& 0.00681\\ \hline
$n^{50}$, Specifically for $n=1$ to 30,000 & 30,000& 3.7 & 0.0048& 0.00256\\ \hline
\end{tabular}
}

\ms \nt{\bf 3.4 \ \ Demographic Statistics and other social
science data} \ss The demographic data that the author used in
this section are obtained from Energy Information Administration
United States Department of Energy. \ss

{\small \cl{Table 7. Analysis of demographic data}
\hspace{0.6in}\begin{tabular}{|l|l|l|l|l|}\hline
& Sample Size &$\chi^2$ &$d_1$ &$d_{\max}$\\ \hline
GDP total for each &1606 &33.0 &0.1009 &0.03352 \\
country in 1992-2001 &&&& \\ \hline
World Consumption of &480 &40.5 &0.1113 &0.05936 \\
Primary Energy by Selected &&&& \\
Country Groups, 1992-2001 &&&& \\ \hline
Population for each &212 &6.59 &0.0792 &0.04160 \\
country in year 2001 &&&& \\ \hline
Crude Oil Production, &4196 &32.2 &0.0383 &0.02166 \\
Import and Export, &&&& \\
Stock Build Statistics for &&&&\\
Each Country in 1990-2001 &&&& \\ \hline
\end{tabular}
}

\ss The chi-square test rejects the null hypothesis (the data's
obeys the Benford's distribution) for all except the population
statistics; while the total variation distance says that the oil
statistics fit the Benford's distribution better than the other
three. \ss Other two sets of statistics in social science obtained
from {\it The World Affair Companion} by Gerald Segal (1996) are
also analyzed. The null hypothesis is accepted in both cases in
the chi-square test. \ss

{\small \cl{Table 8. Analysis of social science data}
\hspace{0.6in}\begin{tabular}{|l|l|l|l|l|}\hline
&Sample Size &$\chi^2$&$d_1$&$d_{\max}$\\ \hline
Annual average number of  &333 &11.05 &0.0715 &0.03222\\
people reported killed or &&&&\\
affected by disasters per &&&&\\
region and country in 1968-1992&&&&\\ \hline
Weapon imports of 50 &291 &11.06 &0.0832 &0.02872\\
leading recipients in&&&&\\
US dollars in 1990-1994 &&&&\\ \hline
\end{tabular}
}

\ms \nt{\bf 3.5 \ \ Numbers appeared in Magazines and Newspaper}
\ss Among the data compiled in Benford's original paper [2],
``numbers appeared in Reader's Digest'' attracts our attention.
Numbers appeared in a magazine or newspaper seems random at first
glance; its conformance to Benford's law was deemed as a
coincidence rather a rule. However, by referring to the Hill's
paper, ``selected at random and random samples $\ldots$ then the
significant-digit frequencies of the combined sample will converge
to the logarithmic distribution'' [11], we hypothesize that
numbers appeared in Reader's Digest could be considered as a
mixture of data of different distributions and more specialized
magazines containing biased mixture of data will deviate from
Benford's distribution more significantly. Here numbers appeared
in both Innovation (science orientated) and The Economist (current
affairs orientated, more on social side) were analyzed. It was
found that neither conforms well to the Benford's distribution;
nonetheless, the combination dataset conforms to Benford much
better with total variation distance of 0.0474. This supports our
hypothesis. \ss

{\small \cl{Table 9. Analysis of numbers appeared in magazines}
\hspace{0.6in}\begin{tabular}{|l|l|l|l|}\hline
Numbers appeared in &Sample Size &$x^2$ &$d_1$  \\ \hline
Innovation (Vol. 2 No. 4. 2002) &152    &6.57   &0.0747\\ \hline
The Economist (17th May 2003)   &449    &16.24 &0.0788 \\ \hline
Combination of the above two    &601    &9.22    &0.0474\\ \hline
\end{tabular}
}

\ms \nt {\bf 3.6 \ \ Other Datasets} \ss Other kinds of datasets,
which are usually quoted as example of Benford's law include:
survival distribution [14], the magnitude of the gradient of an
image \& Laplacian Pyramid Code [15], radioactive half-lives of
unhindered alpha decays [11]. Datasets compiled by Benford in his
original paper include: Rivers area, population, newspapers,
specific heat, molecular weights, atomic weights, Reader's Digest,
X$-$Ray volts, American league, addresses, death rate $\ldots$ [2]
They provide a wide range of examples to illustrate Benford's Law;
however, the certainty of their conformance to Benford's law
varies largely. \bs \nt {\bf 4. \ \ Application of Benford's Law}
\ss The analysis of first few significant digit frequencies
provides us a potential framework to examine the accuracy and
authenticity of data values in numerical data set.

\ss \nt {\bf 4.1 \ \ The Design of Computers and Analysis of
Roundoff Errors} \ss Peter Schatte has determined that based on an
assumption of Benford input, the computer design that minimizes
expected storage space (among all computers with binary-power
base) is base 8 [16], and other researchers have started exploring
the use of logarithmic computers to speed up calculations [17]. It
will be crucial for the development of computing science in the
information age. \ss \nt {\bf 4.2 \ \ A Goodness-of-fit against
Benford to Detect Fraud} \ss A$.$ Rose and J$.$ Rose developed a
VBA code for detecting fraud under Excel working environment [8].
But the program can only do a first digit test but not a first
significant digit test (for example, the 1 in 0.150 will not be
counted by the program), and the program will encounter error upon
reading text (not numbers). We improved the program to eliminate
these limitations of the original program. \ss However, we argue
that the use of a goodness-of-fit test against Benford to detect
fraud is dubious. It has many limitations. Firstly, this method to
detect fraud will not be accurate if the data in the account has
very close built-in maximums and minimums. Secondly, a one-off
embezzlement or very few entries of amendments made in an
accounting record will not be detected using Benford's Law.
Moreover, the public or specifically here the fraudsters learn
very quickly; after a few cases of successful detections of fraud
using Benford's law being reported, they may take special cares
upon changing the entries in the account table in order to conform
to Benford's law. Therefore, we conclude that Benford's law can
only serve as an initial strategy to detect fraud and it will
become less useful later. \ss \nt{\bf  4.3. \ \ Others} \ss Mark
Nigrini mentioned in his paper [4] that Benford's law could be
used to test the accuracy of measuring equipment. Hill in his
paper suggested Benford's law as a test of reasonableness of
forecasts of a proposed model [9]. \ms \nt{\bf 5. \ \ Conclusions}
\ss Among the large range of samples we have analyzed, by using
Chi-square test and calculating total variation distance and
maximum of deviations, Fibonacci series, Sequence $\alpha^n$,
population and number in magazines fit Benford's distribution
well. Fitness to Benford's law of other data such as physical
constants and stock prices is controversial, and it requires
further investigations. Two common explanations for the origin of
Benford's distribution are multiplicative process and scale
invariance. One simple multiplicative process cannot reproduce
most of naturally occurring phenomena thus by virtue of
multiplicative process it cannot explain the origin of Benford's
distribution for all datasets that are fit. The understanding of
the origin of scale invariance has been one of the fundamental
tasks of modern statistical science. How systems with many
interacting degrees of freedom can spontaneously organize into
scale invariant states will be the focus of further research. \ss
It is generally believed that the more chaotic and diversified the
probability distributions are, the better the overall data set
fits the Benford's Law. Thus, we put forward a bold conjecture
that Benford's Law is an act of nature and can be used to indicate
the randomness of our world. Therefore, it can be related to the
second thermodynamics law of entropy. This will be an interesting
area for further investigation. \bs \nt {\bf Acknowledgement} \ss
We would like to thank Assoc Professor Choi Kwok Pui from the
Department of Mathematics, National University of Singapore for
his help and guidance in this project. \bs

\end{document}